\documentclass{article}

\usepackage[english]{babel}

\usepackage[letterpaper,top=2cm,bottom=2cm,left=3cm,right=3cm,marginparwidth=1.75cm]{geometry}

\usepackage{amsmath}
\usepackage{graphicx}
\usepackage[colorlinks=true, allcolors=blue]{hyperref}
\usepackage{amsthm}
\usepackage{booktabs}
\usepackage{enumerate}
\newtheorem{theorem}{Theorem}
\newtheorem{lemma}{Lemma}
\newtheorem{corollary}{Corollary}
\newtheorem{conjecture}{Conjecture}
\newtheorem{question}{Question}

\title{On a Conjecture of Sun Zhi-Wei and Related Diophantine Equations}
\author{Wang Jia-Hui and Zhu Hui-Lin\\
School of Mathematical Sciences, Xiamen University,\\
Xiamen City 361005, Fujian Province, P.R.China\\}

\begin{document}
\maketitle

\begin{abstract}
For any integer $m\ge0$, we recall that triangular numbers are those $\mathbf{T}(m)=\frac{m(m+1)}{2}$. A conjecture of Sun Zhi-Wei states that an integer $2^n\pm n$ with any $n>2$ can not be a triangular number. The motivation of this work is to confirm this conjecture.
\end{abstract}

\section{Introduction}
  In MathOverflow(MO403375) Sun Zhi-Wei asked his question: Is there an integer $n>2$ with $2^{n}-n$ or $2^{n}+n$ a triangular number? I guess the answer should be negative, but I don't know how to prove this. Any ideas? We recall triangular numbers are those
$$T(m)=\frac{m(m+1)}{2}(m=0,1,2,\cdots).$$
Clealy,
$$2^{1}-1=1=T(1),\ 2^{1}+1=3=T(2),\ 2^{2}+2=6=T(3).$$
Also see the remark 8.6 of Conjecture 8.6 in the book\cite{ref1}.
%Sun's conjecture, essentially, is to determine whether Diophantine equation with form $2^n \pm g(n)=\frac{m(m+1)}{2}$ where $g(n)\in \mathbf{Z[n]}$, has solution or not. The equation is equivalent to  $x^2+D(n)=2^{n+3}$ where $x=2m+1$ and $D(n)=\mp 8g(n)-1$. It is a kind of generalized Ramanujan-Nagell equation.
\newline \newline We can transform the question into a generalized Ramanujan-Nagell equation as the following:$$2^{n}\pm n=T(m)=\frac{m(m+1)}{2},$$
$$2^{n+3}\pm 8n=4m(m+1)=(2m+1)^{2}-1,$$
$$(2m+1)^{2}\pm 8n-1=2^{n+3}.$$

Let $x=2m+1$ be odd and $D(n)=\pm 8n-1$. We have

\begin{equation}
\label{eq1}x^{2}+D(n)=2^{n+3},D(n)\equiv7(\bmod8).
\end{equation}
It is a generalized Ramanujan-Nagell equation of $D(n)$ being a integer coefficient polynomial of $n$ more than a constant.
\newline \newline In this work, by applying the deep results of Bauer and Bennett about Hypergeometric method\cite{ref2}, we prove the equation (1) has no solution with $n>2$ when $D(n)=\pm 8n-1$, i.e. we answered Sun's conjecture positively as the following:
\begin{theorem}
Equation
\begin{equation}
    \label{eq1}2^n\pm n=\frac{m(m+1)}{2}
\end{equation}
has no solution with $n\in \mathbf{Z},n>2$ and $m \in \mathbf{Z},m\ge0$.
\end{theorem}
Our method is also effective for some kind of equations of the form  $x^2+D(n)=y^n,$ where $D(n)\in \mathbf{Z}(n),y$ is a given integer more than 1. For example, we prove the following theorem:
\begin{theorem}
Equation
\begin{equation}
    \label{eq2}3^n\pm (n^3+n)=x^2
\end{equation}
has only one positive integer solution $3^1-(1^3+1)=1^2$.
\end{theorem}

\section{Preliminaries}
We mainly used the result of M. Bauer and M. A. Bennett as the following
%Before we proceed our work, we need state Bauer and Bennett's result in approximation of quadratic irrational number.

\begin{lemma}{[2, Corollary 1.6]}

If \(y\) is a positive integer in the table below, then, if $p \in \mathrm{Z}$ and $q=y^k$, for some integer $k>2$, with
$$(y,k)\notin \{(2,3),(2,7),(2,8),(3,7)\}.$$
we have
$$\left|\sqrt{y}-\frac{p}{q}\right|>q^{-\lambda_2(y)}$$
where $\lambda_2(y)$ is as follows:

\begin{table}[h!]
\begin{center}
\caption{}
\begin{tabular}{cccccccccccc}

\toprule
\(y\) & $\lambda_2(y)$ & \(y\) & $\lambda_2(y)$ & \(y\) & $\lambda_2(y)$ & \(y\) & $\lambda_2(y)$ & \(y\) & $\lambda_2(y)$ & \(y\) & $\lambda_2(y)$ \\
\midrule
2&1.48&21&1.67&38&1.72&53&1.51&69&1.73&83&1.55\\
3&1.65&23&1.45&40&1.55&54&1.29&70&1.75&84&1.46\\
5&1.36&24&1.64&42&1.61&55&1.39&72&1.58&85&1.43\\
6&1.46&26&1.97&43&1.91&56&1.76&73&1.51&87&1.89\\
10&1.99&28&1.64&44&1.68&57&1.76&74&1.69&90&1.98\\
12&1.65&29&1.60&45&1.53&58&1.66&75&1.91&91&1.82\\
13&1.53&30&1.91&46&1.21&60&1.47&76&1.27&92&1.58\\
14&1.84&31&1.70&47&1.66&62&1.58&77&1.44&93&1.73\\
17&1.94&33&1.73&48&1.63&63&1.76&78&1.70&95&1.95\\
18&1.87&34&1.74&50&1.81&65&1.76&79&1.56&96&1.41\\
19&1.87&35&1.87&51&1.65&66&1.57&80&1.72&98&1.54\\
20&1.67&37&1.55&52&1.81&68&1.46&82&1.71&99&1.69\\
\bottomrule
\end{tabular}
\end{center}
\end{table}
\end{lemma}
By the hypergeometric method from Pad$\acute{e}$ approximation, Bauer and Benett improved the upper bound of $n$ in the equation (1) from $\frac{10\log|D|}{\log2}+435$\cite{ref3} to $5.55\log\left(|D|\right)$\cite{ref2} for given constant $D$.

\iffalse
Let $y$ be $2$, we have
\begin{corollary}
If $p\in \mathbf{Z}$ and $q=2^k$ for some $k\in \mathbf{Z}$, then
$$\left| \sqrt{2}-\frac{p}{q}\right| >q^{-1.48}$$
\end{corollary}
Our goal is to prove
\begin{theorem}
Equation
\begin{equation}
    \label{eq1}2^n\pm n=\frac{m(m+1)}{2}
\end{equation}

Has no solution with $n\in \mathbf{Z},n>2$ and $m \in \mathbf{Z},m\ge0$.
\end{theorem}

\begin{corollary}
Equation
\begin{equation}
    \label{eq2}2^n\pm n^2=\frac{m(m+1)}{2}
\end{equation}
Has no solution with $n\in \mathbf{Z},n>2$ and $m \in \mathbf{Z},m\ge0$.
\end{corollary}
\fi

\section{Proof of Theorem 1}

%Let  $x$ be $2m+1$ and $D$ be $\mp 8n-1$, then the equation (1) is equivalents to
%\begin{equation}
    %\label{eq3}x^2+D=2^{n+3}
%\end{equation}.

\begin{description}
\item[Case 1:]$D(n)=-8n-1$.\newline
\begin{equation}
\label{eq2}x^2=2^{n+3}+8n+1.
\end{equation}
If \(n\) is odd, then there exists $t \in \mathbf{Z}\ $ such that $n=2t+1$. Hence

$$x^2=2^{2t+4}+16t+9.$$

When $t=1$ or $t=2$, it is not true since both 89 and 297 are not perfect squares.
\newline When $t\ge 3$, from $16t+9< 2^{t+3}+1$ we have
$$(2^{t+2})^2=2^{2t+4}<x^2<2^{2t+4}+2^{t+3}+1=(2^{t+2}+1)^2.$$
So there is no integer solution for equation (4).\newline
If \(n\) is even, suppose $n=2t(t\ge2)$, then

$$x^2=2^{2t+3}+16t+1.$$

When $t=2$ or $t=3$, it is not true since both 161 and 561 are not perfect squares.\newline
When $t\ge 4$, from $2^{t+\frac{5}{2}}+1>16t+1$ we have
$$(2^{t+\frac{3}{2}})^2=2^{2t+3}<x^2<(2^{t+\frac{3}{2}}+1)^2=2^{2t+3}+2^{t+\frac{5}{2}}+1,$$
$$x=[2^{t+1}\sqrt{2}]+1,$$
where $[2^{t+\frac{3}{2}}]$ is the Gauss integer function.
Suppose
$$f(t)=([2^{t+1}\sqrt{2}]+1)^2-2^{2t+3}.$$
We will prove $f(t)\neq 16t+1$ for any $t\ge 4$.
One may verify that when $t\in \{4,5,6,7,8,9,10\}$, $f(t)\neq 16t+1$.\newline
Suppose $t\ge 11$. %we will prove $f(t)>16t+1$ by induction.
%$t=12$, then $f(13)=24329>16\times 13+1=209$.\newline
Naturally, $\left\{2^{t+1}\sqrt{2}\right\}=2^{t+1}\sqrt{2}-[2^{t+1}\sqrt{2}]$, then

$$f(t)=(2^{t+1}\sqrt{2}-\{2^{t+1}\sqrt{2}\}+1)^2-2^{2t+3}$$
$$=\{2^{t+1}\sqrt{2}\}^2+1-2^{t+2}\sqrt{2}\{2^{t+1}\sqrt{2}\}+2^{t+2}\sqrt{2}-2\{2^{t+1}\sqrt{2}\}$$
$$=2^{t+2}\sqrt{2}(1-\{2^{t+1}\sqrt{2}\})+(1-\{2^{t+1}\sqrt{2}\})^2.$$
Consider $$\left|\frac{1-\{2^{t+1}\sqrt{2}\}}{2^{t+1}}\right|$$
$$=\left|\frac{1-2^{t+1}\sqrt{2}+[2^{t+1}\sqrt{2}]}{2^{t+1}}\right|$$
$$=\left|\frac{1+[2^{t+1}\sqrt{2}]}{2^{t+1}}-\sqrt{2}\right|.$$
By Lemma 1,
$$\left|\frac{1+[2^{t+1}\sqrt{2}]}{2^{t+1}}-\sqrt{2}\right|>\frac{1}{2^{1.48(t+1)}}.$$
Hence,
$$2^{t+2}\sqrt{2}(1-\{2^{t+1}\sqrt{2}\})+(1-\{2^{t+1}\sqrt{2}\})^2>2^{t+2}\sqrt{2}(1-\{2^{t+1}\sqrt{2}\})$$
$$=2^{2t+3}\sqrt{2}\left|\frac{1-\{2^{t+1}\sqrt{2}\}}{2^{t+1}}\right|>\frac{2^{2t+3+\frac{1}{2}}}{2^{1.48(t+1)}}=2^{0.52t+2.02}>16t+1(t\ge11).$$
This indicates $f(t)>16t+1$ for any $t\ge 11$.
\item[Case 2:]$D(n)=8n-1$.\newline
If $n$ is odd, similarly with Case 1, then there is no integer solution.% suppose $n=2t+1$, then we get
%$$x^2=2^{2t+4}-16t-7$$
%When $t=1$, $t=2$ or $t=3$, it is not true since 41, 217 or 969 are non-perfect squares.
%When $t\ge4$, we have
%$$2^{t+2}-1<x<2^{t+2}$$
%$x$ can not be an integer, though there is no solution when $n$ is odd.
\newline
If $n$ is even, suppose $n=2t(t\ge2)$, then
$$x^2=2^{2t+3}-16t+1.$$
When $t=2$ or $t=3$, it is not true since both 97 and 465 are not perfect squares.\newline
When $t\ge 4$, we have
$$2^{t+\frac{3}{2}}-1<x<2^{t+\frac{3}{2}}$$ and  $x=[2^{t+\frac{3}{2}}]$.
Now we will prove $x\neq[2^{t+\frac{3}{2}}]$.
\newline
Suppose $f(t)=2^{2t+3}-[2^{t+\frac{3}{2}}]^2$.
Similarly, $f(t)=2^{t+\frac{5}{2}}\left\{2^{t+\frac{3}{2}}\right\}-\left\{2^{t+\frac{3}{2}}\right\}^2$. Calculation validates $f(t)\neq 16t-1$ when $t\in \{4,5,6,7,8,9,10\}$. So\newline
$$f(t)>2^{t+\frac{5}{2}}\left\{2^{t+1}\sqrt{2}\right\}-1.$$
By Lemma 1, we have
$$f(t)>2^{2t+\frac{7}{2}}\left(\sqrt{2}-\frac{\left[2^{t+1}\sqrt{2}\right]}{2^{t+1}}\right)-1>2^{0.52t+2.02}-1>16t-1(t\ge11).$$
In conclusion, equation (1) has no integer solution.
\end{description}

\section{Proof of Theorem 2}
\begin{description}
\item[Case 1:]$x^2=3^n+n^3+n.$\newline
If $2|n$, suppose $n=2t(t\ge2)$, then
$$3^t<x<3^t+1(t\ge7).$$
And for $2\leq t\leq 6$, there is no integer solution for equation (3).\newline
When $2\not|n$, we suppose $n=2t+1(t\ge1)$. If $t\ge29$, then
$$x^2=3^{2t+1}+(2t+1)^3+2t+1,$$
$$3^{t+\frac{1}{2}}<x<3^{t+\frac{1}{2}}+1,$$
$$x=\left[3^{t+\frac{1}{2}}\right]+1,$$
$$x^2-3^{2t+1}=2\cdot3^{t+\frac{1}{2}}\left(1-\left\{3^{t+\frac{1}{2}}\right\}\right)+\left(1-\left\{3^{t+\frac{1}{2}}\right\}\right)^2>2\cdot3^{2t+\frac{1}{2}}\left|\sqrt{3}-\frac{1+\left[3^{t+\frac{1}{2}}\right]}{3^t}\right|.$$
By lemma 1, we have
$$x^2-3^{2t+1}>2\cdot3^{0.35t+0.5}>(2t+1)^3+2t+1.$$
All values of $3^{2t+1}+(2t+1)^3+2t+1(1\leq t\leq28)$ is as follows(Table 2):
\begin{table}[h!]
\begin{center}
\caption{}
\begin{tabular}{cc}

\toprule
    $t$ & $3^{2t+1}+(2t+1)^3+2t+1$ \\
\midrule
	 1&57\\
     2&373\\
     3&2537\\
     4&20421\\
     5&178489\\
     6&1596533\\
     7&14352297\\
     8&129145093\\
     9&1162268345\\
     10&10460362485\\
     11&94143191017\\
     12&847288625093\\
     13&7625597504697\\
     14&68630377389301\\
	 15&617673396313769\\
16&5559060566591493\\
17&50031545099042617\\
18&450283905891048053\\
19&4052555153019035625\\
20&36472996377170855365\\
21&328256967394537157177\\
22&2954312706550833789813\\
23&26588814358957503391657\\
24&239299329230617529707781\\
25&2153693963075557766443449\\
26&19383245667680019896945653\\
27&174449211009120179071336937\\
28&1570042899082081611640719813\\
     \bottomrule
\end{tabular}
\end{center}
\end{table}

No perfect square in the table, so there is no solution of equation (3) under this condition.

\item[Case 2:]$x^2=3^n-n^3-n.$\newline
If $2|n$, similarily with Case 1, there is no solution for equation (3).\newline
If $2\not|n$, let $n=2t+1(t\ge0)$, when $t>28$,
$$3^{t+\frac{1}{2}}-1<x<3^{t+\frac{1}{2}}$$
This shows $x=\left[3^{t+\frac{1}{2}}\right]$, and we have
$$3^{2t+1}-x^2=2\cdot3^{t+\frac{1}{2}}\left\{3^{t+\frac{1}{2}}\right\}-\left\{3^{t+\frac{1}{2}}\right\}^2>2\cdot3^{2t+\frac{1}{2}}\left|\sqrt{3}-\frac{\left[3^{t+\frac{1}{2}}\right]}{3^t}\right|-1.$$
By Lemma 1, we have
$$3^{2t+1}-x^2>2\cdot3^{0.35t+0.5}-1>(2t+1)^3+2t+1.$$ So there is no solution with $t\ge29$. Check there is only one solution with $0\leq t\leq 28$, that is $t=0,n=1$.
\end{description}

\section{Remarks}

\iffalse
As discussed previously, the object of study for this paper is a kind of generalized Ramanujan-Nagell equation
\begin{equation}
\label{eq6}x^2+D=2^n,
\end{equation}
where $D(n)$ is an integer coefficient polynomial of $n$, and such equation have at most finite solutions by hypergeometric method. We present some conjectures as the following:\newline
\begin{conjecture}
Equation (5) has at most 5 positive integer solutions, where $D$ is an integer coefficient polynomial of $n$.
\end{conjecture}
\fi

From the proof of the two theorems we will see, for any given integer coefficient polynomial $D(n)$, with the degree $k$, the generalized Ramanujan-Nagell equation
\begin{equation}
    \label{eq5}x^2+D(n)=2^n,
\end{equation}
has at most finite solutions, which can be proved by the upper bound of $n$ from Pad$\acute{e}$ approximation.\newline
\begin{question}
For given integer coefficient polynomial $D(n)$ with degree $k$, equation (5) has at most $N(k)$ positive integer solutions. For any given $k$, what is the $N(k)$?
\end{question}
For arbitrary integer coefficient polynomial $D(n)$, we can construct some examples, in which there are arbitrary many positive integer solutions of equation (5). Actually, we find a method to construct the polynomials $D(n)$.\newline
Let
$$D_1(n)=(n-3)(n-4)(n-5)+7,$$
by simple computation, we know
$$x^2+D_1(n)=2^n$$
has at least 3 solutions:
$$(x,n)=(1,3),(3,4),(5,5).$$
Set
$$D_2(n)=(n-3)(n-4)(n-5)(n-6)+c_2(n-3)(n-4)(n-5)+7,$$
Let $n=6$. By solving the equation
$$x^2\equiv 7-2^6(\bmod 6),$$
we have $x\equiv3(\bmod 6)$. Replace $x$ by 3 one has $c_2=8$ (actually, $x$ has infinitly many choices, and the same $c_2$), namely
$$D_2(n)=(n-3)(n-4)(n-5)(n-6)+8(n-3)(n-4)(n-5)+7,$$
and the equation
$$x^2+D_2(n)=2^n$$
has at least 4 solutions:
$$(x,n)=(1,3),(3,4),(3,6),(5,5).$$
Similarly, we have
$$D_3(n)=(n-3)(n-4)(n-5)(n-6)(n-7)+c_3(n-3)(n-4)(n-5)(n-6)+8(n-3)(n-4)(n-5)+7,$$
$$D_4(n)=(n-3)(n-4)(n-5)(n-6)(n-7)(n-8)+c_4(n-3)(n-4)(n-5)(n-6)(n-7)$$
$$+c_3(n-3)(n-4)(n-5)(n-6)+8(n-3)(n-4)(n-5)+7,$$
$$D_5(n)=(n-3)(n-4)(n-5)(n-6)(n-7)(n-8)(n-9)+c_5(n-3)(n-4)(n-5)(n-6)(n-7)(n-8)$$
$$+c_4(n-3)(n-4)(n-5)(n-6)(n-7)+c_3(n-3)(n-4)(n-5)(n-6)+8(n-3)(n-4)(n-5)+7.$$
We may choose $c_3=-3, c_4=1, c_5=-6$, and get the equations
$$x^2+D_3(n)=2^n,$$
$$x^2+D_4(n)=2^n,$$
$$x^2+D_5(n)=2^n,$$
have at least 5, 6 and 7 solutions separately, which are respondingly the following:
$$(x,n)=(1,3),(1,7),(3,4),(3,6),(5,5),$$
$$(x,n)=(1,3),(1,7),(3,4),(3,6),(3,8),(5,5),$$
$$(x,n)=(1,3),(1,7),(3,4),(3,6),(3,8),(5,5),(65,9).$$
\newline
It is worth noting that the values of $c_2, c_3, c_4$ and $c_5$ are not unique. Further more, they satisfy some congruence conditions.

\iffalse
We present the following conjecture:
\begin{conjecture}
For given positive integer $k$ and any integer coefficient polynomial $D$ with $deg(D)\leq k$, the supremum of the number of solutions of equation (5) exists , denoted as $N(k)$. Besides $N(k)\ge k$.
\end{conjecture}

And we may ask:
\begin{question}
Is there a functional relationship between $k$ and $N(k)$? If yes, then what is it?
\end{question}
\fi

Specially, when $D(n)$ is a polynomial with degree 0, i.e. a constant with $D(n)\equiv7(\bmod 8)$, many authors have already given some results of the number of solutions(let $N(D)$ denote the number of solutions):
\begin{theorem}{[3,Theorem 2]}
Let $D$ be a positive integer of $D\equiv7(\bmod 8)$. The equation (5)
has at most one solution except
%two or more solutions in positive integers $x,n$ if and only if $D=7,23$ or $2^k-1$ for some $k\ge4$. The solutions in these exceptional cases are given by
\begin{description}
\item[(I)$D(n)=7,$]$$(x,n)=(1,3),(3,4),(5,5),(11,7),(181,15),$$
\item[(II)$D(n)=23,$]$$(x,n)=(3,5),(45,11),$$
\item[(III)$D(n)=2^k-1(k\ge4),$]$$(x,n)=(1,k),(2^k-1,2k-2).$$
\end{description}
\end{theorem}
\begin{theorem}{[4,Theorem 1 and Theorem 2]}
Let $D$ be a negative integer of $D\equiv7(\bmod 8)$. The equation (5)
has at most 3 solutions except
%\begin{equation}
%x^2-D=2^{n+2},x>0,n>0,D\in \mathbf{N},2\not|D.
%\end{equation}
%And three types of $D$:
\begin{description}
\item[(I)]$$D=-(2^{2m}-3\cdot2^{m+1}+1),\ m\ge3,$$
$$(x,n)=(2^{m}-3,3),(2^{m}-1,m+2),(2^{m}+1,m+3),(3\cdot2^{m}-1,2m+3).$$
\end{description}
Specially, it is found there are 2 cases in which the equation has exactly 3 integer solutions as the following:
\begin{description}
\item[(II)]$$D=-\left(\frac{2^{2m-1}-17}{3}\right)^{2}+32,\ m\ge4,$$
$$(x,n)=\left(\frac{2^{2m-1}-17}{3},5\right),\left(\frac{2^{2m-1}+1}{3},2m+1\right),\left(\frac{17\cdot2^{2m-1}-1}{3},4m+3\right).$$
\item[(III)]$$D=-(2^{2m_{1}}+2^{2m_{2}}-2^{m_{1}+m_{2}}-2^{m_{1}+1}-2^{m_{2}+1}+1),\ \ m_{1}>m_{2}+1>2,$$
$$(x,n)=(2^{m_{1}}-2^{m_{2}}-1,m_{2}+2),(2^{m_{1}}-2^{m_{2}}+1,m_{1}+2),(2^{m_{1}}+2^{m_{2}}-1,m_{1}+m_{2}+2).$$
\item[(IV)] And if $D$ is not one of the above types and the equation $$u'^2-Dv'^2=-1$$ has solutions $(u',v')$, then $N(D)\leq2$.
\end{description}
\end{theorem}
We give the following questions:
\begin{question}
What is the necessary and sufficient condition that equation (5) has no solution including two cases of $D$ being positive and negative? In which condition, equation (5) has exactly 1, 2, 3 solutions when $D$ is negative?
\end{question}
Actually some authors researched more general Ramanujan-Nagell equation as the following form:
\begin{equation}
x^2+D=k\cdot2^n,
\end{equation}
where $k$ is a fixed odd integer and $D$ is a given constant (maybe not $D\equiv7(\bmod8)$). For example, J. Stiller\cite{ref5} find the equation
$$x^2+119=15\cdot2^n$$
has exactly six solutions, which are
$$(x,n)\in \{(1,3),(11,4),(19,5),(29,6),(61,8),(701,15)\}.$$\
M. Ulas\cite{ref6} find another two equations
$$x^2-117440512=57\cdot2^n,$$
and
$$x^2-26404=165\cdot2^n,$$
have exactly six solutions, which are
$$(x,n)\in \{(10837,0),(10880,14),(11008,16),(13312,20),(32768,24),(45056,25)\},$$\
and
$$(x,n)\in \{(163,0),(178,5),(218,7),(262,8),(442,10),(838,12)\}.$$\
separately.
So we present another question:
\begin{question}
How many integer solutions the equation (6) has at most?
\end{question}

\subsection*{Acknowledgements}
The second author was supported by China National Nature Foundation Grant(No. 11501477), the Science Fund of Fujian Province(No. 2015J01024), and the Fundamental Research Funds for the Central University (No. 2072017001). The authors thank professor Sun Zhi-Wei for providing this conjecture. The second author thanks professor Cao Wei for inviting him to give the talk in "2021 Zhangzhou Diophantine Equations Conference". The authors thank Chen Han, Hong Hao-Jie and Liu Chang for their help in giving some comments.

\end{document}